\input amstex
\documentstyle{amsppt}
\magnification1200
\pagewidth{6.5 true in}
\pageheight{9.25 true in}
\NoBlackBoxes

\topmatter
\title Multiplicative functions in arithmetic progressions
\endtitle
\author Antal Balog, Andrew Granville and K. Soundararajan
\endauthor
\address{Alfr\'ed R\'enyi Institute of Mathematics, POB 127, 1364 Budapest,
Hungary}\endaddress
\email{balog{\@}renyi.hu}
\endemail
\address{D{\'e}partment  de Math{\'e}matiques et Statistique,
Universit{\'e} de Montr{\'e}al, CP 6128 succ Centre-Ville,
Montr{\'e}al, QC  H3C 3J7, Canada}\endaddress
\email{andrew{\@}dms.umontreal.ca}
\endemail
\address{Department of Mathematics, 450 Serra Mall, Bldg. 380, Stanford University,
Palo Alto CA 94305, USA} \endaddress \email{ksound{\@}math.stanford.edu}
\endemail
\thanks{The first author would like to acknowledge the support of the
Hungarian National Science Foundation, as well as of the
Centre de Recherches en math\'ematiques in Montreal.
Le deuxieme auteur est partiellement soutenu par une bourse
de la Conseil de recherches en sciences naturelles et en g\' enie
du Canada. The third  author is partially supported by the
National Science Foundation and the American Institute of
Mathematics (AIM).}
\endthanks
\abstract We develop a theory of multiplicative functions (with
values inside or on  the unit circle) in arithmetic progressions
analogous to the well-known theory of primes in arithmetic
progressions.
\endabstract
\endtopmatter

\parindent=20pt
\document
\def\cbar{\overline{\chi}}
\def\sumstar{\sideset\and ^* \to \sum}
\def\L{\fracwithdelims()}

 \head 1.
Introduction \endhead

\noindent A central focus of multiplicative number theory has been
the question of how primes are distributed in arithmetic
progressions (see e.g. [3]):  There are $\pi(x)$ primes up to $x$,
of which $\pi(x;q,a)$ belong to the arithmetic progression $a \pmod
q$, and
$$
\pi(x;q,a) \sim \frac{\pi(x)}{\phi(q)}  \tag{1.1}
$$
as $x\to \infty$, whenever $(a,q)=1$. In other words the primes are
eventually equidistributed amongst the plausible arithmetic
progressions $\pmod q$. In applications it is  important to know how
big $x$ must be, as a function of $q$, for (1.1) to hold.
Unconditionally we can only prove (1.1) for all $q$ when $x$ is
enormous, larger than exponential of a (small) power of $q$, whereas
we believe it holds when $x$ is just a little bigger than $q$, that
is $x>q^{1+\epsilon}$.  We can prove a far better range for $x$ if
we ask whether (1.1), or a weakened version of it, holds for most
$q$:\ Fix $\epsilon>0$. There exists $A$ such that $$ \left|
\pi(x;q,a) -  \frac{\pi(x)}{\phi(q)}  \right|
 \leq
 \epsilon \  \frac{\pi(x)}{\phi(q)} \tag{1.2}
$$
for all $(a,q)=1$ for all $q\leq x^{1/A}$, except {\sl possibly}
those $q$ that are multiples of some {\sl exceptional} modulus $r$
(which depends only on $\epsilon$ and $x$).  We do not believe that
any such $r$ exists but if it does then $r \geq \log x$ so the
result (1.2) applies to most $q$. Moreover if $q$ is divisible by
$r$ we can still get an understanding of the distribution of primes
in the arithmetic progressions mod $q$, albeit rather different from
what we had expected:\ There exists a real character $\psi \pmod
{r}$ such that
$$
\left( \pi(x;q,a)-\frac{\pi(x)}{\phi(q)} \right) = \psi(a) \left(
\pi(x;q,1)-\frac{\pi(x)}{\phi(q)} \right) + O\left( \epsilon
\frac{\pi(x)}{\phi(q)}\right),  \tag{1.3}
$$
whenever $(a,q)=1$ and $r$ divides $q$, with $q\leq x^{1/A}$.

In the literature (e.g. [3, 15]) exceptional  moduli are derived and
discussed in terms of $L$-functions and egregious counterexamples to
the Generalized Riemann Hypothesis known as Siegel-Landau zeros. Our
goal here is to develop this same theory without $L$-functions and
in a way that generalizes to a wide class of functions.  The first
step of our approach comes in recognizing that results about the
distribution of prime numbers can usually be rephrased as results
about mean values of the M\"obius function $\mu(n)$ or the Liouville
function  $\lambda(n)$, both multiplicative functions (see, e.g.
[19]), as we run through certain sequences of integers $n$. To see
this write the Von Mangoldt function   as
$$
\Lambda (n) = \sum_{d|n} \mu(d) \log (n/d)
$$
(where $\Lambda (p^e)=\log p$ for prime $p$ and integer $e\geq 1$,
and $\Lambda(n)=0$ otherwise), which is deduced from $\log =
1*\Lambda$ by  M\"obius inversion.  Now  $\log$ is a very ``smooth''
function, so that estimates for the sum of $\Lambda (n)$ for $n$ in
an arithmetic progression, follow from estimates for the mean value
of $\mu$ in certain related arithmetic progressions. Thus (1.1) is
equivalent to the statement that the mean value of $\mu$ in any
arithmetic progression $\pmod q$ tends to 0 as we go further and
further out in the arithmetic progression [19]; and indeed (1.2) is
tantamount to the fact that
$$
\left| \sum\Sb n\leq x \\ n\equiv a \pmod q\endSb \mu(n) \right|
\leq \epsilon \frac xq \tag{1.2b}
$$
for all $(a,q)=1$ for all $q\leq x^{1/A}$, except  those $q$ that
are multiples of some   exceptional  modulus $r$.

If $q$ {\sl is} a multiple of $r$ then the corresponding result
for $\mu(n)$ is more enlightening than (1.3).  It is well-known that
if (1.3) holds then $L(1,\psi)$ is surprisingly small, which can
hold only if $\psi(p)=-1$ for many of the ``small'' primes $p$.  In
other words $\mu(p)=\psi(p)$ for many small primes $p$, which
implies that $\mu(n)$ displays a bias towards looking like
$\psi(n)$. But then $\mu(n)$ tends to look like $\psi(n)=\psi(a)$
for $n\equiv a {\pmod q}$, so its mean value over such integers $n$
tends to be like $\psi(a)$. All of these vagaries can be made more
precise, but for now we content ourselves with the appropriate
analogy to (1.3):
$$
 \sum\Sb n\leq x \\ n\equiv a \pmod q\endSb \mu(n)
 =\psi(a) \sum\Sb n\leq x \\ n\equiv 1 \pmod q\endSb \mu(n)
  +O\left( \epsilon \frac xq \right)  \tag{1.3b}
$$
whenever $(a,q)=1$ and $r$ divides $q$, with $q\leq x^{1/A}$. The
key point is to see that $\mu$ {\sl pretends} to be $\psi$, at least
at small values of the argument in this exceptional case. This
notion generalizes very well.

This paper is concerned with determining, for multiplicative
functions $f$ with $|f(n)|\le 1$  for all $n$, estimates for the
mean-values
$$
\frac{1}{x/q} \ \sum\Sb n\le x\\ n\equiv a\pmod q\endSb f(n)
\tag{1.4}
$$
when $(a,q)=1$.
We will show that for any fixed $\epsilon>0$ there exists $A$ such that
$$
\left| \sum\Sb n\leq x \\ n\equiv a \pmod q\endSb f(n) \right| \leq
\epsilon \frac xq \tag{1.2c}
$$
for all $(a,q)=1$ for all $q\leq x^{1/A}$, except  possibly  those
$q$ that are multiples of some  exceptional  modulus $r$.
Moreover, if such a modulus $r$ exists there is a character $\psi
\pmod {r}$ such that if $q$ is a multiple of $r$  with $q\leq
x^{1/A}$, then
$$
 \sum\Sb n\leq x \\ n\equiv a \pmod q\endSb f(n)
 =\psi(a) \sum\Sb n\leq x \\ n\equiv 1 \pmod q\endSb f(n)
 +O\left( \epsilon \frac xq \right)  \tag{1.3c}
$$
whenever $(a,q)=1$.

There are examples of $f$ for which exceptional characters {\sl do}
exist: indeed we simply let $f(n)=\psi(n)$ or even
$f(n)=\psi(n)n^{it}$ for some $|t|\ll 1/\epsilon$. Thus, in this
theory, exceptional characters take on a different role in that they
exist for any $f$ which pretends to be a function of the form
$\psi(n)n^{it}$. This pretentiousness can be made more precise in
terms of the following {\sl distance function} $\Bbb D=\Bbb D_1$
between two multiplicative functions $f$ and $g$ with $|f(n)|,
|g(n)|\leq 1$: \footnote{This distance function is particularly
useful since it satisfies the triangle inequality  ${\Bbb
D}_q(f_1(n),g_1(n);x) + {\Bbb D}_r(f_2(n),g_2(n);x) \ge {\Bbb
D}_{qr}((f_1f_2)(n),(g_1g_2)(n);x)$, which may be proved by squaring both
sides and using the Cauchy-Schwarz inequality (see Lemma 3.1 of
[12], or the general discussion in [13]).}
$$
{\Bbb D}_r(f(n),g(n); x)^2 :=\sum\Sb p\le x \\ p\nmid r\endSb
\frac{1-\text{Re } f(p)\overline{ g(p)}}{p}.
$$

\subhead{1.2.  More precise results}\endsubhead

Throughout $f$ is a multiplicative function with $|f(n)|\le 1$ for
all $n$, and $x\geq Q,A\geq 1$ are given. Of all primitive
characters with conductor below $Q$, let $\psi \pmod r$ be that
character  for which
$$
\min_{|t| \le A} \ {\Bbb D}_r(f(n),  \psi(n) n^{it} ;x)^2 \tag{1.5}
$$
is a minimum.\footnote{If there are several possibilities for
$\psi$, simply pick one of those choices.}. Let $t=t(x,Q,A)$ denote
a value of  $t$ that gives the minimum value in (1.5).

\proclaim{Theorem 1}  With the notation as above we have, when $a
\pmod q$ is an arithmetic progression with $(a,q)=1$ and $q\leq Q$,
that
$$
\sum\Sb n\le x\\ n\equiv a\pmod q \endSb f(n) \ll \text{\rm
Error}(x,q,A)   \tag{1.2d}
$$
if $r$ does not divide $q$, and
$$
\sum\Sb n\le x\\ n\equiv a\pmod q \endSb f(n) =
\frac{\psi(a)}{\phi(q)} \sum\Sb n\le x \\ (n,q)=1\endSb f(n)
\overline{\psi(n)} + O\Big( \text{\rm Error}(x,q,A)\Big) \tag{1.3d}
$$
if $r$ does divide $q$. Here we may take either $Q=x^{\frac 1A}$
with $\log x \ge A\ge 20$ and
$$
\text{\rm Error}(x,q,A)= \frac{x}{q} \cdot  \frac{1}{\sqrt{\log A}}
;
$$
or $Q=\log x$ and $A=\log^2 x$ with
$$
\text{\rm Error}(x,q,\log^2 x)=\frac{x}{q (\log x)^{1/3+o(1)}} +
\frac{x}{(\log x)^{1-o(1)}} .
$$
\endproclaim

Throughout let $\chi_0$ be the primitive character mod $q$; and here
let $\chi=\chi_0$ if $r$ does not divide $q$, and  $\chi=\psi\chi_0$
if $r$ does divide $q$. Then (1.2d) and (1.3d) can be expressed in
one equation as
$$
 \sum\Sb n\leq x \\ n\equiv a \pmod q\endSb f(n) \ - \ \chi(a) \sum\Sb n\leq x \\
 n\equiv 1 \pmod q\endSb f(n) \ll   \text{\rm Error}(x,q,A) . \tag{1.3e}
$$
We can deduce the following useful bound from Corollary 2.2 below,
which is useful in bounding (1.3d) of Theorem 1:
$$
\frac{1}{\phi(q)} \sum\Sb n\le x \\ (n,q)=1\endSb f(n) \overline{\psi(n)}
\ll \frac{x}{q} \left( (1+{\Bbb D}_q(f(n),  \psi(n) n^{it} ;x)^2)
e^{-{\Bbb D}_q(f(n),  \psi(n) n^{it} ;x)^2}  + \frac 1{(\log x)^{1/4}}
\right). \tag{1.6}
$$

We see that (1.2c) and (1.3c) follow  immediately from this Theorem.
The theorem gives more than (1.2c) and (1.3c) as we get results here
with some uniformity, and somewhat  stronger results when $x$ is
larger than $e^q$.

We suspect that  the constant ``$1/3$'' in
Theorem 1 can be increased, but we can
show that it cannot be taken to be any larger than ``$1$'': Assume that
the Generalized Riemann Hypothesis holds in order to make the
calculations easier. Suppose that $f(p)=\psi(p)$ for all primes
$p\leq y$ and $f(p)=\chi(p)$ for primes $p\in (y,x]$ for some
distinct characters $\psi\pmod r$ and $\chi \pmod \ell$ with
$(r,\ell)=1$. If $y=x/2$ then (1.4) equals
$\psi(a)(x-(\pi(x)-\pi(x/2)))/\phi(q)
 + O(\sqrt{x}\log^2x)$
if $r|q$, it equals $\chi(a)(\pi(x)-\pi(x/2))/\phi(q) + O(\sqrt{x}\log^2x)$
if $\ell|q$, and $O(\sqrt{x}\log^2x)$ otherwise; so that the  upper bound
given on Error$(x,q,A)$
 cannot be improved beyond $x/(\phi(q)\log x)$.
If $y=x^{o(1)}$ and $\log y\geq  (\log x)^{1-o(1)}$ then one can show that the
two main terms in the calculation above are
$\sim \psi(a)\Psi(x,y)/\phi(q)$ and
$\sim \chi(a)\Phi(x,y)L(1,\chi_0\psi\overline\chi)/\phi(q)$
where $\Psi(x,y)$, and $\Phi(x,y)$, denote the number of integers up to $x$
free of prime factors $>y$, and $\leq y$, respectively. Selecting
$\log y=2\log x \log\log\log x/\log\log x$ we deduce that
the  upper bound
given on Error$(x,q,A)$
 cannot be improved beyond
$x\log\log x/(\phi(q)\log x\log\log\log x)$.

Earlier work in this area   focussed on the equidistribution of $f(n)$ in
arithmetic progressions; that is,
on obtaining bounds for
$$
\sum\Sb n\le x \\ n\equiv a\pmod q \endSb f(n) - \frac{1}{\phi(q)} \sum\Sb n\le x\\ (n,q) =1
\endSb f(n)   \tag{1.7}
$$
whenever $(a,q)=1$. Elliott [4,5] showed that this is  $\ll x (\log
\log x/\log x)^{\frac 18}$ for all $q$ except possibly for multiples
of a certain exceptional modulus $r$: note that this is non-trivial
in the range $q\le (\log x)^{\frac 18+ o(1)}$. Hildebrand [14]
showed that (1.7) is $\ll x/(q\sqrt{\log A})$ for all $q\le Q$
except possibly for multiples of one of at most two exceptional
moduli $r$ and $r'$. Our  result improves both of these, by
understanding the asymptotics in all cases, at worst in terms of the
exceptional moduli.

One can see many of the ideas in this paper, and particularly the
dichotomy between the two cases, appearing in work on elementary
proofs of the prime number theorem. In particular Selberg [18]
essentially shows that there is no exceptional modulus  by giving
(what is equivalent to) an elementary proof that the minimum in
(1.5) cannot be too small. This is developed in [8] more along the
lines given here.

\medskip

One might also ask whether
$$
\lim_{x\to \infty} \ \frac{1}{x/q} \sum\Sb n\le x \\ n\equiv a\pmod
q \endSb f(n) \tag{1.8}
$$
exists, for a given progression $a\pmod q$ with $(a,q)=1$? If
$f(n)=n^{it}$ for some $t\ne 0$, or $\psi(n)n^{it}$ for some
character $\psi \pmod q$ then the answer is evidently ``no''; or
even if $f(n)$ pretends to be $\psi(n)n^{it}$.

\proclaim{Corollary 1}   If $f$ is unpretentious, that is ${\Bbb
D}(f(n), \psi(n)n^{it} ;\infty)=\infty$ for all characters $\psi$
and real numbers $t$ then
$$
\sum\Sb n\le x\\ n\equiv a\pmod q \endSb f(n) = o(x)  \tag{1.9}
$$
as $x\to \infty$, for any $a\pmod q$ with $(a,q)=1$.   If there does
exist a primitive character $\psi \pmod r$ and a real number $t$ for
which  ${\Bbb D}(f(n), \psi(n)n^{it} ;\infty)<\infty$ then $\psi$
and $t$ are unique. In this case (1.9) still holds if $r$ does not
divide $q$. However, if $r$ divides $q$ then
$$
\sum\Sb n\le x\\ n\equiv a\pmod q \endSb f(n) = (1+o(1))
\frac{\psi(a) x^{1+it}}{q(1+it)} \prod_{p \nmid q}
\Big(1-\frac{1}{p}\Big) \Big( 1+
\frac{f(p)\overline{\psi(p)}}{p^{1+it}} +
\frac{f(p^2)\overline{\psi(p^2)}}{p^{2+2it}} +\ldots\Big).
$$
as $x\to \infty$.
In particular if $f$ is real-valued and $\chi$
is that real character $\pmod q$ for which
$$
\prod_{p}\Big(1-\frac 1p \Big) \Big( 1+\frac{f(p)\chi(p)}{p} +\frac{f(p^2)\chi(p^2)}{p^2}
+ \ldots \Big)
$$
is  maximized,\footnote{Note that $|\prod_{p\leq x} (1-1/p)(1+g(p)/p+g(p^2)/p^2+ \dots)|
\asymp \exp (-{\Bbb D}(1, g(n),x)^2 )$.}  then
$$
\lim_{x \to \infty} \ \frac{1}{x/q} \sum\Sb n\le x\\ n\equiv a\pmod
q \endSb f(n) = \chi(a) \prod_{p\nmid q} \Big(1-\frac{1}{p}\Big)
\Big(1+\frac{f(p)\chi(p)}{p} +\frac{f(p^2)\chi(p^2)}{p^2}+
\ldots\Big).
$$
\endproclaim

\subhead{1.3. Mean values of multiplicative functions, twisted by Dirichlet characters}\endsubhead

Traditionally one estimates mean values of functions in an arithmetic progression
by using characters, as in the identity
$$
\sum\Sb n\le x\\ n\equiv a\pmod q \endSb f(n) = \frac{1}{\phi(q)}
\sum_{\chi \pmod q} \chi(a) \sum\Sb n\le x \endSb f(n) \cbar(n).
\tag{1.10}
$$
We therefore give estimates on such characters sums.

\proclaim{Theorem 2}  Define $\psi$ as above.
For fixed  $Q \le \sqrt{x}$ and $A=\log^2x$ the estimate
$$
\sum\Sb n\le x\endSb f(n) \cbar(n) \ll x \Big(\frac{\log Q}{\log x} \Big)^{\frac 1{20}},
$$
holds for all characters $\chi$ of conductor $q\le Q$, except perhaps
those induced from $\psi$.
Moreover the estimate
$$
 \sum_{n\le x} f(n) \overline{\chi}(n) \ll \frac{x}{(\log x)^{1/3+o(1)}}.
$$
holds for all characters $\chi$ of conductor $q\le \log x$, except perhaps
those induced from $\psi$.
\endproclaim

When $q$ is small compared to $x$ (that is $q\le \log x$)
we may obtain
good estimates on (1.4) by using Theorem 2 in (1.10).
Moreover Theorem 2 for large moduli $q$ suggests that  Theorem 1 should
hold for these moduli, but deriving this from (1.10) is not
straightforward since many characters are now involved.

We now record one more consequence of this circle of ideas.  Let $f$
be a  real-valued completely multiplicative function with $|f(n)|\le
1$ for all $n$. Proving a conjecture of Hall, Heath-Brown, and
Montgomery, Granville and Soundararajan showed in [10] that
$$
\sum_{n\le x} f(n) \ge (\delta_1 +o(1)) x,
$$
uniformly for all $f$, where
$o(1)\to 0$ as $x\to \infty$, and
$$
\delta_1 = 1 -2\log (1+\sqrt{e}) + 4\int_1^{\sqrt{e}} \frac{\log t}{t+1} dt = -0.656999\ldots .
$$
Further this lower bound is best possible, and is attained when $f(\ell)=1$
for
all primes
$\ell\le x^{1/(1+\sqrt{e})}$, and $f(\ell)=-1$ for $x^{1/(1+\sqrt{e})}\le \ell\le x$.
For any totally multiplicative function $g$ with each $g(\ell)=1$ or $-1$
and any $x$, there exists infinitely many primes $p$ for which
$\L {\ell}{p}=g(\ell)$ for all primes $\ell\leq x$, as may be proved using
quadratic reciprocity and Dirichlet's theorem for primes in arithmetic
progression. Thus taking $g=f$ we see that  there exist primes $p$
such that
at least $17.15 \%$ of the integers below $x$ are quadratic residues $\pmod p$.
Here $17.15 \%$ is an approximation to $\delta_0 = (1+\delta_1)/2 = 0.1715\ldots$,
and this constant is best possible.  We now give a generalization of this
result to arithmetic progressions.

\proclaim{Corollary 2}   Let $a \pmod q$ be a progression with
$(a,q)=1$.  Then
$$
\liminf_{x\to \infty}\  \inf\Sb p \\ p\nmid q\endSb \
 \frac{1}{x/q}
 \sum\Sb n\le x\\ n\equiv a\pmod q\endSb \L{n}{p} =
 \cases \delta_1 &\text{if  } a \equiv \square \pmod q\\
 -1 &\text{if } a \not\equiv \square \pmod q.\\
 \endcases
 $$
 Colloquially, if $a \equiv \square \pmod q$ then at least $17.15\%$ of the integers
 $n\le x$ with $n \equiv a\pmod q$ are quadratic residues $\pmod p$ for
 any prime $p$.
 \endproclaim

{\bf Acknowledgments.}   We are grateful to Ben Green for some
interesting conversations regarding Proposition 5.3, and for drawing
our attention to [1], and to John Friedlander for drawing our
attention to [7].

\head 2. Hallowed Hal\'asz \endhead

Our main results rest  on the large sieve and development of
Hal\'asz's pioneering results on mean values of multiplicative
functions, given here after incorporating significant refinements
due to Montgomery and Tenenbaum [20]  (and see [11] for an explicit
version).

\proclaim{Theorem 2.1 (Hal\'asz)}  For any $T\ge 1$ we have
$$
\frac 1x \ \sum_{n\le x} f(n) \ll   \ \max_{|t| \le T} \
\Big(1+{\Bbb D} (f(n),n^{it};x)^2\Big)
e^{-{\Bbb D}(f(n),n^{it};x)^2}  \
+ \frac{1}{\sqrt{T}}.
$$
\endproclaim

\noindent Hal\'asz's result shows that if the mean-value of $f$ is
large then  $f(p)$ pretends to be $p^{it}$ for some small value of
$t$, the quantity ${\Bbb D}$ giving an appropriate measure of the
distance between $f(p)$ and the values $p^{it}$.  Note that if
$f(p)=p^{it}$ then the mean-value of $f$ is indeed large since
$\sum_{n\le x} n^{it}  \sim x^{1+it}/(1+it)$.

We need the following consequence of Theorem 2.1:

\proclaim{Corollary  2.2}  For $1\leq T\leq (\log x)^{\frac 12}$
select $t$ as in Theorem 2.1. If
$r$ is an integer $\leq \sqrt x$ then
$$
\frac 1{\frac{\phi(r)}r\ x} \ \sum\Sb n\le x \\ (n,r)=1 \endSb  f(n)
\ll   \Big(1+{\Bbb D}_r (f(n),n^{it};x)^2\Big)
e^{-{\Bbb D}_r(f(n),n^{it};x)^2}  \ + \frac{1}{\sqrt{T}}.
$$
\endproclaim

Hal\'asz's theorem allows us to estimate $\sum_{n\le x} f(n)
\cbar(n)$, and deduce that this mean-value is small unless $f$
pretends to be $\chi(p)p^{it}$ for some small $t$.  We will show in
Lemma 3.1 below that   $f$ can pretend to be ${\chi(p)}p^{it}$ for
at most one character $\chi \pmod q$ provided $q$ is not too large, and
hence  Hal\'asz's theorem  implies Theorem 2 for large moduli $q$.

Proving an old conjecture of Erd{\H o}s and Wintner, Wirsing [21]
showed that every real valued multiplicative function $f$ with
$|f(n)|\le 1$ has a mean-value; and it is natural to ask whether
this holds as $n$ ranges over an arithmetic progression: that is,
whether (1.8) exists?  We determined this mean value in Corollary 1
above, and then gave the analogy to Wirsing's result  for arithmetic
progressions. These results  rest on a further beautiful result of
Hal\'asz, which gives  a qualitative version of Theorem 2.1.\footnote{
Although unrelated to the main topic of this paper, we take this
opportunity to mention a beautiful, but not widely known consequence
of these results of Hal\'asz and Wirsing. Let $G$ be a finite abelian
group, and let $f: {\Bbb N} \to G$ be multiplicative.  Then for any
$g\in G$ the density $\lim_{x\to \infty} \frac 1x \#\{n\le x: f(n)
=g\}$ exists!  This is stated as a Conjecture ``of the sixteen year
old Hungarian mathematician I. Ruzsa" by Erd{\H o}s [6].  }

\proclaim{Theorem 2.3 (Hal\'asz)}  If there exists real number $t$ for
which ${\Bbb D} (f(n),n^{it};\infty)<\infty$ then
$$
\sum_{n\le x} f(n) = (1+o(1))\frac{x^{1+it}}{1+it} \prod_{p} \Big( 1-\frac 1p\Big) \Big(1+\frac{f(p)p^{-it}}{p}
+ \frac{f(p^2)p^{-2it}}{p^2} +\ldots\Big),
$$
as $x\to \infty$,
where the $o(1)\to 0$ in a manner that depends on $f$.
If ${\Bbb D} (f(n),n^{it};\infty)=\infty$ for all $t$ then $\sum_{n\le x} f(n)=o(x)$ as $x\to \infty$
\endproclaim

Let $S$ be a subset of the unit disc such that  $|\text{arg}(1-z)|
\le \theta < \pi/2$ for all $z\in S$ (here we take the argument to
always be between $-\pi$ and $\pi$). We can also deduce from Theorem
2.3 that if $f$ is such that $f(p)\in S$ for all $p$, then the limit
in (1.8) exists (extending Theorem 2).

\demo{Proof of Corollary 2.2} By Theorem 4  of [11] and the last two
displayed equations in the proof of Corollary 3 of [11],\footnote{We
take this opportunity to correct an error in the statement of those
two equations: in each, one should divide the sum on the right side
of the equation by the number of terms in the sum.} we have, for $t$
selected with ${\Bbb D}_r (f(n),n^{it};x)$ minimal and $d\leq
\sqrt{x}$
$$
\sum_{n\leq x/d} f(n) = \frac 1{d^{1+it}} \sum_{n\leq x} f(n)
+O\left( \frac xd \ \frac{\log\log x}{(\log x)^{2-\sqrt{3}}} \right)
.
$$
We deduce, from the combinatorial sieve  that, for $r\leq \sqrt{x}$,
$$
\frac{r}{\phi(r)}   \sum\Sb n\leq x \\ (n,r)=1 \endSb f(n) =
\prod_{p|r} \left( 1-\frac{f(p)}{p^{1+it}}\right)  \left(
1-\frac{1}{p}\right)^{-1} \sum_{n\leq x} f(n) + O\left(\left(
\frac{r}{\phi(r)}\right)^2 x \ \frac{\log\log x}{(\log
x)^{2-\sqrt{3}}} \right) .
$$
There is no loss of generality in setting $f(p)=p^{it}$ for each
prime $p$ dividing $r$, and so we deduce the result since
$2-\sqrt{3}>1/4$.
\enddemo

\head 3.  Pretentious characters are repulsive \endhead

We show that $f$ cannot pretend to be two different functions of the
form $\psi(n)n^{it}$:

\proclaim{Lemma 3.1} For each  primitive character $\psi$ with
conductor below $\log x$ select $|t| \le \log^2 x$ for which ${\Bbb
D}(f(n),\psi(n)n^{it};x)$ is minimal, and then label these pairs so
that  $(\psi_j, t_j)$ is that pair which gives the $j$-th smallest
distance ${\Bbb D}(f(n),\psi_j(n)n^{it_j};x)$. Let $r_j$ be the conductor
of $\psi_j$. For each $j\geq 1$ we
have
$$
 {\Bbb D}_{r_j}(f(n),\psi_j(n)n^{it_j};x)^2 \ge \Big(1-\frac{1}{\sqrt{j}} \Big) \log \log x +O(\sqrt{\log \log x}).
$$
\endproclaim

\demo{Proof}  This is essentially Lemma 3.4 of [13], except that
there we dealt with the case where $f$ itself was a character, and
took all $t_j=0$. The same proof applies, and for completeness we
give the details.  Note that
$$
\align {\Bbb D}(f(n),\psi_j(n)n^{it_j};x)^2 &\ge \frac{1}{j} \sum_{k
=1}^j {\Bbb D}(f(n),\psi_k(n)n^{it_k};x)^2 \ge \frac{1}{j} \sum_{p\le
x} \frac{1}{p } \sum_{k=1}^{j} (1-\text{Re
}f(p)\overline{\psi_k(p)p^{it_k}} ) \\
& \ge \frac{1}{j} \sum_{p\le x} \frac 1p
\Big(j - \Big|\sum_{k=1}^{j } \psi_k(p)p^{it_k} \Big|\Big).\\
\endalign
$$
By Cauchy-Schwarz we have that
$$
\Big( \sum_{p\le x} \frac 1p \Big| \sum_{k=1}^{j}
\psi_k(p)p^{it_k}\Big|\Big)^2 \le \Big( \sum_{p\le x} \frac 1p\Big)
\Big( \sum_{p\le x}  \frac 1p  \Big| \sum_{k=1}^{j}
\psi_k(p)p^{it_k}\Big|^2 \Big).
$$
The first factor above is $= \log \log x +O(1)$, while the second term is
$$
\le \sum_{p\le x} \frac 1p \Big( j+ \sum\Sb 1\le k,\ell \le j \\ k\neq \ell \endSb
\psi_k(p)\overline{\psi_\ell(p)}p^{i(t_k-t_\ell)} \Big) = j \log \log x +O(j^2),
$$
upon using the prime number theorem in arithmetic progressions. For
the transition from $\Bbb D$ to $\Bbb D_r$ simply note that
$\sum_{p|r} 1/p \ll \log\log\log\log x$. The Lemma follows.
\enddemo

\proclaim{Corollary  3.2} If $m$ is a positive integer for which
$f^m$ is real valued and ${\Bbb D}(f(n),\psi(n)n^{it};x)^2 \le
\frac{1}{16m^2} \log \log x$ with $|t|\leq \frac 1m \log^2x$ then
$\psi^m$ is real and $|t|\ll \frac 1{m\sqrt{\log x}}$.
\endproclaim

\demo{Proof} Let $m=1$. Taking complex conjugates we have  ${\Bbb
D}(f(n),\psi(n)n^{it};x)={\Bbb D}(f(n),\overline\psi(n)n^{-it};x)$. By
Lemma 3.1 this is impossible unless $\overline\psi=\psi$, that is
$\psi$ is real. Now, by the triangle inequality we have ${\Bbb
D}(f(n)^2,n^{2it};x)\leq 2{\Bbb D}(f(n),\psi(n)n^{it};x)$. As $f^2\geq 0$
we have ${\Bbb D}(f(n)^2,n^{2it};x)^2\geq \sum_p 1/p$ where the sum is
over those primes $p\leq x$ with $\cos(2t\log p)<0$. Since we get
the ``expected'' number of primes in intervals $[y,y+y^{1-\delta}]$
for some fixed $\delta>0$ we can deduce that $\sum_p 1/p \geq \frac
12 (\log\log x -\log (\max\{ 1/|t|, \log |t|\}  ) +O(1))$ for
$|t|\log x\geq 1$, which implies that $|t|\ll 1/\sqrt{\log x}$, as
required. For larger $m$ we use the triangle inequality to note that
${\Bbb D}(f(n)^m,\psi^m(n)n^{imt};x)\leq m{\Bbb D}(f(n),\psi(n)n^{it};x)$
and the result then follows from the $m=1$ case.
\enddemo

\proclaim{Lemma 3.3} With the hypothesis of Lemma 3.1 we have
$$
{\Bbb D}_{r_2}(f(n),\psi_2(n)n^{it_2};x)^2 \ge \left\{ \frac{1}{3}
+o(1)\right\} \log \log x .
$$
\endproclaim

\demo{Proof} The proof of Lemma 3.1 gives
$$
{\Bbb D}(f(n),\psi_2(n)n^{it_2};x)^2  \ge  \sum_{p\le x} \frac 1p
\Big(1 -  \frac{1}{2} \Big|\sum_{k=1}^{2 } \psi_k(p)p^{it_k}
\Big|\Big) = \sum_{p\le x} \frac 1p \Big(1 -  \frac{1}{2}
\Big|1+\chi(p)p^{it}  \Big|\Big) ,
$$
taking $\chi=\psi_2\overline{\psi_1}$ and $t=t_2-t_1$. By the prime
number theorem for arithmetic progressions, this equals
$$
\frac 1{\phi(q)} \sum\Sb a \pmod q \\ (a,q)=1 \endSb \ \int_1^{\log
x}   \Big(1 -  \frac{1}{2} \Big|1+\chi(a)e^{itv}\Big| \Big) \ \frac
{\text{d}v}v + o(\log\log x) .
$$
If $\chi$ has order $m>1$ then there are exactly $\phi(q)/m$ values
of $j \pmod m$ for which $\chi(a)=e^{2i\pi j/m}$, and so our
integral equals, taking $v=2u$,
$$
 \frac 1m \sum_{j=0}^{m-1}   \ \int_1^{\log x} \Big(1 -
\Big| \cos( tu + \pi j/m) \Big| \Big) \ \frac {\text{d}u}u
+o(\log\log x) . \tag{3.1}
$$

Now consider $\int_1^{\log x} | \cos( tu + \beta) | du/u$ for
arbitrary $\beta$. If $u\leq \epsilon /t$ then $\cos( tu + \beta)
=\cos( \beta) +O(\epsilon)$. The contribution of $u$ in the range $
\epsilon /t<u<1/(\epsilon t)$ to the integral is $O(1)$. For $u\geq
1/(\epsilon t)$ we cut the range up into intervals $[U, U+2\pi/t)$,
and use the fact that $1/u=1/U+O(1/U^2t)$, to obtain
$$
\int_U^{U+2\pi/t} | \cos( tu + \beta) | \ \frac {\text{d}u}u =
\int_U^{U+2\pi/t} \left\{  \frac 1{2\pi} \int_{w=0}^{2\pi} |\cos w|
dw  + O(1/tu) \right\} \frac {\text{d}u}u .
$$
The total error that arises  here is $\ll \int_{u\geq 1/(\epsilon
t)} du/(tu^2) \ll \epsilon$. We deduce that (3.1) equals $\log\log
x$ times
$$
\alpha \cdot \frac 1{2\pi} \int_{w=0}^{2\pi} (1-|\cos w|) dw +
(1-\alpha) \cdot  \frac 1m \sum_{j=0}^{m-1} (1 - |\cos(\pi j/m)|)
+o(1)
$$
for some $\alpha$ in the range $0\leq \alpha\leq 1$. Now $\frac
1{2\pi} \int_{w=0}^{2\pi} (1-|\cos w|) dw = 1-\int_{-1/2}^{1/2}
\cos( \pi \theta) d \theta=1-2/\pi$, and $1-\frac 1m
\sum_{-m/2<j\leq m/2} \cos( \pi j/m)=1-1/(m \sin(\pi/2m)$ if $m$ is
odd, $1-1/(m \tan(\pi/2m)$ if $m$ is even. The minimum occurs for
$m=3$, and the result follows.
\enddemo

Our next lemma (which is essentially Proposition 7 of [13])
formulates a similar repulsion principle involving characters of
substantially larger conductors, at the cost of obtaining a smaller
separation.

\proclaim{Lemma 3.4} Let $\chi \pmod q$ be a non-principal character
(with $q\ge 3$), and $t\in \Bbb R$. There is an absolute  constant
$c>0$ such that for all $x\ge q$ we have
$$
 {\Bbb D}_q(1,\chi(n)n^{it};x)^2 \ge
 \frac12 \log \Big( \frac{c\log x}{\log (q(1+|t|))} \Big) .
$$
Consequently, if $f$
is a multiplicative function, and $\chi$ and $\psi$ are any two
characters with conductors $q,r\leq Q$, respectively, such that
$\chi\overline{\psi}$ is non-principal, then for $x \ge Q$ we have
$$
{\Bbb D}_q(f(n),\chi(n)n^{it};x)^2 + {\Bbb D}_r(f(n),\psi(n) n^{iu};x)^2
\ge \frac{1}{8} \log \left( \frac{c\log x}{2\log (Q(1+|t-u|))}
\right) .
$$
\endproclaim
\demo{Proof}  For completeness we sketch the proof.  We consider
 $d_{\chi,t}(n) =\sum_{ab=n} \chi(a)a^{it} \cbar(b)b^{-it}$.
The proof of the  P{\' o}lya-Vinogradov inequality is easily
modified to show that $\sum_{n\leq N} \chi(n) \ll (\phi(q)/q)
\sqrt{q}\log q$. Using this and partial summation, we see that
$$
 \sum\Sb n\le x\endSb \chi(n)n^{it} = x^{it} \sum_{n\le x} \chi(n) - it \int_1^x u^{it-1} \sum_{n\le u}\chi(n) du
 \ll \frac{\phi(q)}q \ \sqrt{q} (\log q)(1+|t|\log x).
$$
Therefore, grouping the terms $ab=n$ according to whether $a \le \sqrt{x}$ or $b\le \sqrt{x}<a$, we obtain
$$
\sum_{n\le x} d_{\chi,t}(n)\ll \left( \frac{\phi(q)}q \right)^2 \
\sqrt{qx} (\log q) (1+|t|\log x) . \tag{3.2}
$$

We write $d(n) = \sum_{\ell |n}d_{\chi,t}(n/\ell)h(\ell)$ where $h$
is a multiplicative function with $h(p)=2-2\text{Re
}(\chi(p)p^{it})$, and $|h(n)|\le d_4(n)$ for all $n$.   Then
$$
 x\log x +O(x) = \sum_{n\le x} d(n) = \sum_{\ell \le x} h(\ell) \sum_{n\le x/\ell} d_{\chi,t}(n).
$$
When $\ell \le x/(q^2 (1+|t|)^2)$ we use (3.2) to bound the sum over
$n$. For larger $\ell$ we use that $|d_{\chi,t}(n)|\leq \chi_0(n)
d(n)$ where $\chi_0(n)$ is the principal character mod $q$. Let
$B=(\log q)^{q/\phi(q)}$. Thus if $x/B \leq \ell \leq x$ then
$$
\sum_{n\leq x/\ell} |d_{\chi,t}(n)|\leq \sum_{n\leq x/\ell} d(n) \ll
\frac x \ell \ \log (e x/\ell) \ll \frac x \ell \ \log B \ll \frac x
\ell \ \log q\ (\phi(q) /q)^2,
$$
since $q/\phi(q)\ll \log \log q$. Now, if $N>B$ then
$$
\sum_{n\leq N} \chi_0(n) d(n) = \sum\Sb ab\leq N \\ (ab,q)=1\endSb 1
\ll \sum\Sb a\leq \sqrt{N} \\ (a,q)=1\endSb \sum\Sb b\leq N/a
\\ (b,q)=1\endSb 1 \ll \sum\Sb a\leq \sqrt{N} \\ (a,q)=1 \endSb
\frac{\phi(q)}q \frac Na \ll \left( \frac{\phi(q)}q \right)^2  N
\log N
$$
by the sieve which we use for $\ell$ in the range $x/(q(1+|t|))^2\le
\ell \le x/B$ with $N=x/\ell$. Combining these estimates we obtain
 $$
 \align
 x\log x &\ll \sum_{\ell \le x/(q(1+|t|))^2} |h(\ell)|  \left( \frac{\phi(q)}q \right)^2 \sqrt{qx/\ell} \  (\log q)(1+|t|\log (x/\ell)) \\
&\qquad \qquad + \sum_{x/(q(1+|t|))^2\le \ell \le x} |h(\ell)|
\left( \frac{\phi(q)}q \right)^2 \frac{x}{\ell} \log (q(1+|t|))
 \\
 &\ll  \left( \frac{\phi(q)}q \right)^2 x \log (q(1+|t|)) \sum_{\ell \le x} \frac{|h(\ell)|}{\ell}.
 \\
 \endalign
 $$
 Since $\sum_{\ell \le x} |h(\ell)|/\ell \ll \exp(\sum_{p\le x} |h(p)|/p)
 =\exp(2{\Bbb D}(1,\chi(n)n^{it};x)^2)$ we obtain the first statement of the Lemma.

To obtain the second estimate note that the triangle inequality
gives
 $$
 \align
({\Bbb D}_q(f(n),\chi(n)n^{it};x)+ {\Bbb D}_r(f(n),\psi(n)n^{iu};x))^2
&\ge ({\Bbb D}_{qr}(f(n),\chi(n)n^{it};x)+ {\Bbb
D}_{qr}(f(n),\psi(n)n^{iu};x))^2\\
&\ge \sum\Sb p\le x \\ p\nmid qr\endSb \frac{1-\text{Re }|f(p)|^2
\chi(p)\overline{\psi(p)} p^{i(t-u)}}{p}
\\
&\ge \frac{1}{2} \sum\Sb p\le x \\ p\nmid qr\endSb  \frac{1-\text{Re
}\chi(p)\overline{\psi}(p)p^{i(t-u)}}{p},
\\
\endalign
$$
and now we appeal to the first part of the lemma.
\enddemo

\head 4.  The results for small moduli \endhead

\demo{Proof of Theorem 2 for small moduli} Select $\psi$ as at the
start of section 1.2 with $Q=\log x$ and $A=\log^2x$.  Let
${\chi}^{\prime}$ be the primitive character that induces $\chi
\pmod q$.  By assumption $\chi^{\prime}\ne \psi$.  By Lemma  3.3 we
know that for all $|t| \le \log^2 x$
$$
{\Bbb D}(f(n),\chi^{\prime}(n)n^{it};x)^2 \ge \Big( \frac{1}{3}+o(1)
\Big) \log \log x .
$$
Plainly
$$
{\Bbb D}(f(n),\chi(n)n^{it};x)^2 = {\Bbb
D}(f(n),\chi^{\prime}(n)n^{it};x)^2 + O\Big(\sum_{p|q} \frac 1p\Big)
\ge \Big( \frac{1}{3}+o(1) \Big) \log \log x .
$$
Theorem 2 for small moduli now follows from Hal\'asz's Theorem 2.1.
\enddemo

\demo{Proof of Theorem 1 for small moduli}  From our discussion
above we know that for any $\chi \pmod q$ not induced by $\psi$ we
have (for all $|t|\le \log^2 x$)
$$
{\Bbb D}(f(n),\chi(n)n^{it};x)^2 \ge
 \Big( \frac{1}{3}+o(1) \Big) \log \log x .
$$
By a similar argument, Lemma 3.1 implies that for all but at most
$\sqrt{\log \log x}$ characters $\mod q$ we have
$$
{\Bbb D}(f(n),\chi(n)n^{it};x)^2 \ge \log \log x+O(\sqrt{\log\log x}) ,
$$
for all $|t| \le \log^2 x$. Therefore, from (1.10), Hal\'asz's Theorem 2.1 and
these estimates, it follows that, for $r \nmid q$, we have
$$
 \sum\Sb n\le x \\ n\equiv a\pmod q \endSb f(n)
  \ll \frac{x}{q (\log x)^{1/3+o(1)}} + \frac{x}{(\log x)^{1-o(1)}}.
 \tag{4.1}
 $$
 When $r|q$ there is an extra contribution to (1.10) from the character ${\tilde \psi} \pmod q$
 which is induced by $\psi \pmod r$.  Thus in this case we obtain
 $$
 \sum\Sb n\le x \\ n\equiv a\pmod q \endSb f(n)
 =\frac{\psi(a)}{\phi(q)} \sum_{n\le x} f(n) \overline{{\tilde \psi}(n)}
 + O\Big( \frac{x}{q (\log x)^{1/3+o(1)}}
 + \frac{x}{(\log x)^{1-o(1)}}\Big) ,
 \tag{4.2}
 $$
which gives the result.
  \enddemo

\demo{Proof of Corollary 1} If $f$ is unpretentious then the result
follows from Theorem 1  for small moduli and (1.6); and
similarly if $f$ is pretentious but $r$ does not divide $q$. If
$f$ is pretentious and $r$ divides $q$ then the result follows from
(1.3d) and Theorem 2.3. To see that $\psi$ and $t$ are unique suppose
we have another pair  $\psi'$ and $t'$ so that
$D(\psi'(p)p^{it'},\psi(p)p^{it};\infty)<\infty$ by the triangle inequality,
and thus $L(s,\psi'\overline{\psi})$ has a pole at $1+i(t'-t)$ which is false.
Finally, if $f$ is real-valued then, by taking $x$
arbitrary large in Corollary 3.2,  we discover that $\psi$ is real,
and $t$ is zero, and the result follows.
\enddemo

\demo{Proof of Corollary 2} If $|\sum_{n\le x,\ n\equiv a\pmod q}
\L{n}{p}|\gg x/q$, then  $|\sum_{n\le x,\ n\equiv b\pmod q}
\L{n}{p}|\gg x/q$ whenever $(b,q)=1$ by two applications of (1.3e).
Select $c\pmod q$ so that $\psi(c)$ is as close as possible to $i$
and let $b\equiv a/c \pmod q$; we then deduce from (1.3e) that
$\psi$ must be real since both sums are real (see also Corollary
3.2). Therefore (1.3d) tells us that our sum is $(\psi(a)/\phi(q))
\sum_{n\le x} \L{n}{p} \overline{\psi(n)} +o(x/q)$. This is
$(\psi(a)/q) \sum_{n\le x} g(n) +o(x/q)$ where $g$ is the
multiplicative function with $g(p)= \L{n}{p} \overline{\psi(p)}$ if
$p\nmid q$, and $g(p)=1$ if $p|q$, by Proposition 4.4 of [10] (with
$\phi=\pi/2$ and $\epsilon =\log q/\log x$).

Now suppose that $a \not \equiv \square \pmod q$, so that there is a
real character $\psi \pmod q$ for which $\psi(a)=-1$.  There exist
infinitely many primes $p$ for which $\L{\ell}{p} = \psi(\ell)$ for
any prime $\ell \nmid q$ with $\ell \le x$, by quadratic reciprocity
and Dirichlet's theorem. In this case
$$
   \sum\Sb n\le x\\ n\equiv a\pmod q \endSb \L{n}{p}
   = \psi(a) \Big( \frac{x}{q} +O(1)\Big) = - \frac{x}{q}+ O(1),
$$
and so the result follows.

Now suppose that $a \equiv \square \pmod q$. Then $\psi(a)=1$ for
all real characters $\psi \pmod q$ and so our sum is $(1/q)
\sum_{n\le x} g(n) +o(x/q)$. The result then follows from Corollary
1 of [10], the minimal value being obtained by selecting prime $p$
so that $\L{\ell}{p}=1$ for all primes $\ell \le x^{1/(1+\sqrt{e})}$
and $\L{\ell}{p}=-1$ for  $x^{1/(1+\sqrt{e}) }\le \ell \le x$ as
described in section 1.3, and with $\psi$ being the trivial
character $\pmod 1$, so that
$$
   \sum\Sb n\le x \\ n\equiv a\pmod q \endSb \L{n}{p}
   = \frac{1}{\phi(q)} \sum\Sb n\le x\\ (n,q)=1\endSb \L{n}{p} +o(x)
   = (\delta_1 +o(1)) \frac{x}{q} .
$$
(The reader can either directly verify the last equality, or consult
[10].)
\enddemo

\head 5.  The case of large moduli \endhead

\demo{Proof of Theorem 2 for large moduli}   We will use Theorem 2.1
with $T= (\log x/\log Q)^{\frac 1{10}}$. Suppose that there exists a
character $\psi$ with conductor $\leq Q$ and
$$
\Big| \sum\Sb n\le x \endSb f(n) \overline{\psi(n)} \Big| \ge x \Big( \frac{\log Q}{\log x}\Big)^{\frac 1{20}}.
$$
From Theorem 2.1 we conclude that for some real number $|t|\le T$ we
have
$$
{\Bbb D}(f(n),\psi(n)n^{-it};x)^2 \le \Big(\frac 1{20}+o(1)\Big) \log \frac{\log x}{\log Q}.
$$
From Lemma 3.4 it follows that for any  character $\chi$ with
conductor below $Q$ and $\chi\overline\psi$ non-principal, we have
$$
{\Bbb D} (f(n),\chi(n)n^{-iu},x)^2 \ge \frac{1}{14} \log \frac{\log x}{\log Q},
$$
for all $|u|\le T$.  Appealing to Theorem 2.1, we obtain Theorem 2
for large moduli.

\enddemo

We now embark on the proof of Theorem 1 for large moduli.  It is
convenient to define
$$
 F(x;q,a) =\sum\Sb n\le x \\ n\equiv a\pmod q \endSb  f(n).
$$
We will study $F(x;q,a)$ using the large sieve.  First we show that for
many moduli $r$, the distribution of $f(n)$ for $n \le x$, $n\equiv a\pmod q$,
is uniform in sub-progressions $\pmod r$.

\proclaim{Proposition 5.1} Let $1\le a\le q$. For $\eta >1/\sqrt{\log x}$, the set ${\Cal B}(a,\eta)$ of
bad moduli $1< r\le \sqrt{x/q}$ for which
$$
 \sumstar\Sb \psi\pmod r \endSb \Big| \sum\Sb n\le x /q \endSb f(nq+a) \psi (n) \Big|
 \ge \eta \frac{x}{q},
$$
(where $\sumstar$ denotes, as usual, a sum over primitive characters) satisfies
$$
\sum\Sb r\in {\Cal B}(a,\eta) \endSb \frac{1}{\phi(r)} \leq
\frac{2}{\eta^2}. \tag{5.1}
$$
Let ${\Cal G}(a,\eta)$ denote the set of good moduli $r\le
\sqrt{x/q}$, that is those $r$ which are square-free and not
divisible by any element in ${\Cal B}(a,\eta)$.  If $r \in {\Cal
G}(a,\eta)$ with
 $(r,q)=1$ then
$$
F(x;q,a) = rf(r) F(x/r;q,a/r) + O\Big( \frac{x}{q} \Big(\eta d(r) +
\Big( 1-\frac{\phi(r)}{r}\Big)\Big). \tag{5.2}
$$
Here $a/r$ denotes the residue class $ar^{-1} \pmod q$.
\endproclaim

\demo{Proof} By Cauchy's inequality and then the large sieve (as in
[15], Theorem 7.13) we have
$$
\align \sum\Sb r\le \sqrt{x/q} \endSb\frac{1}{\phi(r)} \left( \
\sumstar_{\psi \pmod r}  \Big| \sum_{n\le x/q}
 f(nq+a)
 \psi(n)\Big|
\right)^2  &\le \sum\Sb r\le \sqrt{x/q} \endSb \ \
\sumstar_{\psi \pmod r}  \Big| \sum_{n\le x/q} f(nq+a)
 \psi(n)\Big|^2 \\
& \leq 2 \frac{x}{q} \sum_{n\le x/q} |f(nq+a)|^2 \leq 2 \Big(\frac xq\Big)^2.
\endalign
$$
The estimate (5.1) follows at once.

Now suppose that $r$ is a good modulus, and consider any progression $b\pmod r$ with $(b,r)=1$.
Since $r$ is square-free we check easily that
$$
\sum\Sb \ell |r \endSb \ \ \sumstar_{\psi \pmod \ell } \overline{\psi(b)} \psi(n)
= \cases
\phi(r/d) &\text{ if  } (n,r)= d, \   \text{and } n\equiv b\pmod {r/d}
\ \text{for some} \ d|r;\\
0 &\text{otherwise}.\\
\endcases
$$
Therefore
$$
\align
\sum\Sb n\le x/q \\ n\equiv b\pmod r \endSb f(nq+a)
&=\frac{1}{\phi(r)}  \sum\Sb \ell | r\endSb \ \  \sumstar_{\psi \pmod \ell} \overline{\psi(b)} \sum\Sb n\le x/q \endSb
f(nq+a)\psi(n)\\
&\hskip .5 in + O\Big( \sum\Sb d|r \\ d> 1\endSb \frac{\phi(r/d)}{\phi(r)} \sum\Sb n\le x/q \\ d|n \\ n\equiv
b\pmod {r/d} \endSb 1\Big).\\
\endalign
$$
The error term above is plainly
$$
\ll \sum\Sb d|r\\ d> 1 \endSb \frac{\phi(r/d)}{\phi(r)} \frac{x}{qr}
\ll \frac{(r-\phi(r))}{\phi(r)} \frac{x}{qr}.
$$
Further, since $r$ is good, if $\ell >1$ and $\ell$ divides $r$ then
$\ell$ is good so that
$$
\sumstar_{\psi \pmod \ell} \overline{\psi(b)} \sum\Sb n\le x/q
\endSb f(nq+a)\psi(n) \ll \sumstar_{\psi \pmod \ell} \Big| \sum\Sb
n\le x/q\endSb f(nq+a) \psi(n) \Big| \le \eta \frac{x}{q}.
$$
We conclude that
$$
\sum\Sb n\le x/q \\ n\equiv b\pmod r \endSb f(nq+a) =
\frac{1}{\phi(r)} \sum_{n\le x/q} f(nq+a) + O\Big(\frac{d(r)}{\phi(r)} \eta \frac xq\Big)
+ O\Big( \frac{x}{q\phi(r)} \Big(1-\frac{\phi(r)}{r}\Big)\Big).
$$

If $(r,q)=1$ then we may choose $b$ so that $bq+a \equiv 0 \pmod r$,
and so the left hand side of the last displayed equation  equals
 $$
 \sum\Sb n\le x \\ n \equiv a \pmod q \\ r|n \endSb f(n)
 = f(r) F(x/r;q,a/r) + O\Big( \frac{x}{qr}\Big(1-\frac{\phi(r)}{r}\Big)\Big).
 $$
Hence, for good moduli $r$ that are coprime to $q$ we have
 $$
 \align
 F(x;q,a) &= f(r)\phi(r) F(x/r;q,a/r) + O\Big(  \frac{x}{q} \Big(\eta d(r)+\Big(1-\frac{\phi(r)}{r}\Big)\Big)\Big)
 \\
& =rf(r) F(x/r;q,a/r) + O\Big(  \frac{x}{q} \Big(\eta d(r)+\Big(1-\frac{\phi(r)}{r}\Big)\Big)\Big),  \\
 \endalign
 $$
proving our Proposition.
\enddemo

\proclaim{Proposition 5.2}  Let $x\ge q^{20}$ and set $A=\log x/\log
q$. Let $a$ and $b$ be any two given reduced residues $\pmod q$.
Then we have
$$
F(x;q,ab) F(x;q,1) = F(x;q,a) F(x;q,b) + O\Big(
\frac{x}{q} \ \frac{1}{\sqrt{\log A}}  \ \max\Sb c\pmod q \\
(c,q)=1\endSb |F(x;q,c)|\Big).
$$
\endproclaim
\demo{Proof}   We may assume that $A$ is large, and  set $\eta=C/
\sqrt{\log A}$ for some sufficiently large constant $C>0$.
 We will produce two numbers $r$ and $s$ with the following properties:

(i) both $r$ and $s$ have at most two prime factors and these prime
factors are $\ge 1/\eta$,

(ii) $r$ and $s$ are both in ${\Cal G}(1, \eta)$,
 ${\Cal G}(a, \eta)$, ${\Cal G}(b, \eta)$ and ${\Cal G}(ab, \eta)$,

(iii) the ratio
 $r/s$ lies between $1-\eta$ and $1+\eta$, and

(iv) $r\equiv bs \pmod q$.
\smallskip

Assuming that this can be done, let us now prove Proposition 5.2.
Four applications of Proposition 5.1 give:
$$
\align
 F(x;q,ab) &= r f(r) F(x/r;q,ab/r) + O(\eta x/q), \tag{5.3} \\
 F(x;q,1) &= sf(s) F(x/s;q,1/s)+ O(\eta x/q), \tag{5.4}\\
 F(x;q,a) &= sf(s) F(x/s;q,a/s)+ O(\eta x/q), \tag{5.5}\\
 F(x;q,b) &= rf(r) F(x/r;q,b/r) + O(\eta x/q).  \tag{5.6}\\
\endalign
$$
By assumption we have $a/s\equiv ab/r \pmod q$, and moreover
$|r/s-1| \le \eta$.  Hence
 $$
 \align
 F(x/s;q,a/s) &= F(x/s;q,ab/r) = F(x/r;q,ab/r) + O(|x/r-x/s|/q)
 \\
 &= F(x/r;q,ab/r) +O(\eta x/(qr)),
 \\
 \endalign
 $$
 and, taking $a=1$,
 $$
 F(x/r;q,b/r) = F(x/s;q,1/s)  +O(\eta x/(qr)).
 $$
 Therefore multiplying (5.3) and (5.4), and (5.5) and (5.6) we obtain that
 $$
 \align
 F(x;q,ab) F(x;q,1) &= rsf(r)f(s) F(x/r;q,ab/r)F(x/s;q,1/s)
 + O\Big(\eta \frac xq \max\Sb c\pmod q\\ (c,q)=1\endSb |F(x;q,c)| \Big)
 \\
 &= F(x;q,a)F(x;q,b) +O\Big(\eta \frac xq \max\Sb c\pmod q\\ (c,q)=1\endSb |F(x;q,c)| \Big)  ,
 \\
 \endalign
 $$
which proves the Proposition.

It remains to show the existence of $r$ and $s$. If there is no
Siegel-Landau zero for a character $\pmod q$ then we can take $r$
and $s$ to be primes. To see this note that if there is no
Siegel-Landau zero then minor modifications to the proof of
Proposition 18.5 of [15] imply that there exist constants $c,B>0$
such that
$$
\sum\Sb z<p\leq (1+\eta)z \\ p\equiv a \pmod q\endSb \log p = \frac
{\eta z}{\phi(q)} \left\{ 1 +O\left(  \frac 1{z^{c/\log q}} + \frac
1{\log z}\right) \right\}  \tag{5.7}
$$
for all $z\geq q^B$ and $(a,q)=1$. We may assume that $A\geq B^2$.
Let ${\Cal P}_1$ denote the set of primes in the interval
$[q^B,\sqrt{x/q}]$ which are in the set ${\Cal G}(1,\eta) \cap {\Cal
G}(a,\eta) \cap {\Cal G}(b,\eta) \cap {\Cal G}(ab,\eta)$. By
Proposition 5.1,
$$
\sum_{p \in {\Cal P}_1} \frac{1}{p} = \sum\Sb q^B \le p \le
\sqrt{x/q}\endSb \frac{1}{p} + O\Big(\frac{1}{\eta^2}\Big).
\tag{5.8}
$$
We divide the interval $[q^B,\sqrt{x/q}]$ into intervals of the form
$(z,(1+\eta)z]$. Any two primes $r$ and $s$ in ${\Cal P}_1 \cap
(z,(1+\eta)z]$ meet criteria (i), (ii) and (iii) above, so it only
remains to satisfy criterion (iv). If criterion (iv) does not hold
then the primes in ${\Cal P}_1 \cap (z,(1+\eta)z]$ lie in at most
$\phi(q)/2$ reduced residue classes $\pmod q$. Therefore, by (5.7)
we must have
$$
\sum\Sb z<p\leq (1+\eta)z \\ p\in {\Cal P}_1 \endSb \log p \leq
\frac {\eta z}{2} \left\{ 1 +O\left(  \frac 1{z^{c/\log q}} + \frac
1{\log z}\right) \right\},
$$
so that
$$
\sum_{p \in {\Cal P}_1} \frac{1}{p} \leq \frac 12 \sum\Sb q^B \le p
\le \sqrt{x/q}\endSb \frac{1}{p} + O(1),
$$
which contradicts (5.8), since $\sum_{q^B \le p \le \sqrt{x/q}} \
1/p \asymp \log A$.

Now assume that there is a real character $\chi \pmod q$ with a
Siegel-Landau zero $\beta$.   Then Proposition 18.5 of [15] yields
that
$$
\sum\Sb  p\leq z \\
p\equiv a \pmod q\endSb \log p = \frac {1}{\phi(q)} \left\{ z -
\chi(a) \frac {z^\beta}{\beta} +O\left(  \frac z{\log^2 z}\right)
\right\} \tag{5.9}
$$
for all $z\geq q^B$ and $(a,q)=1$. Notice that the primes are
concentrated in the residue classes $c \pmod q$ with $\chi(c)=-1$,
and it is therefore difficult to solve $r\equiv bs \pmod q$ in
primes $r$ and $s$ if $\chi(b)=-1$. However products of two primes,
will now be concentrated on the residue classes $c \pmod q$ for
which $\chi(c)=1$, as one may deduce directly from (5.9):
$$
\sum\Sb  p_1p_2\leq z \\ p_1p_2\equiv a \pmod q \\ p_1 > p_2 \geq
q^B\endSb \log p_1 \log p_2 = \frac {1}{\phi(q)} \
\log(\sqrt{z}/q^B) \left\{ z + \chi(a) \frac {z^\beta}{\beta}
+O\left( \frac z{\log z}\right) \right\} , \tag{5.10}
$$
for all $z\geq q^{3B}$ and $(a,q)=1$.

So now let ${\Cal P}_2$ denote the set of products of two primes
$p_1 > p_2 \geq q^B$ with $p_1p_2\leq \sqrt{x/q}$, which are in the
set ${\Cal G}(1,\eta) \cap {\Cal G}(a,\eta) \cap {\Cal G}(b,\eta)
\cap {\Cal G}(ab,\eta)$, and imitate as best we can the proof above.
Proposition 5.1 tells us, again, that ${\Cal P}_2$ contains lots of
elements, and that we can find an interval $(z,(1+\eta)z]$
containing many elements from ${\Cal P}_1$ and ${\Cal P}_2$.  Then
criterion (iv) is met by choosing either two elements from ${\Cal
P}_1$ (when $\chi(b)=1$), or an element each from ${\Cal P}_1$ and
${\Cal P}_2$ (when $\chi(b)=-1$), and the proof of the Proposition
is complete. Although we will not go into the details, the easiest
way to formulate this proof is to combine (5.9) and (5.10) into the
equation
$$
\sum\Sb  p\leq z \\
p\equiv a \pmod q\endSb \log p + \frac 1{\log(\sqrt{z}/q^B)} \ \sum\Sb  p_1p_2\leq z \\
p_1p_2\equiv a \pmod q \\ p_1 > p_2 \geq q^B\endSb \log p_1 \log p_2
= \frac {2z}{\phi(q)}   \left\{ 1  +O\left( \frac 1{\log z}\right)
\right\} , \tag{5.11}
$$
for all $z\geq q^{3B}$ and $(a,q)=1$, removing the effect of the
Siegel-Landau zero, so that we do not need to split this proof into
cases depending on the value of $\chi(b)$.
\enddemo

\remark{Remark}  It is amusing to note that we cannot prove (5.11)
under the assumption that there is no Siegel-Landau zero. In this
case there is an additional $1/z^{c/\log q}$ in the error term, as
in (5.7). It would be interesting to know whether this can be
removed, so as to prove (5.11) uniformly. Note that (5.11) is a
complicated formulation of  Selberg's formula [18] for $z$
sufficiently large.

The proof of Proposition 5.2 is the only proof in this paper that
requires non-trivial information about the zeros of $L$-functions
and it is pleasing to give an alternate proof that is
``elementary''. Thus, in the spirit of Selberg's paper, one can ask
whether one can give an elementary proof of (5.11) in this range?
The proof of Proposition 5.2 goes through easily with  $\sim
2z/\phi(q)$ on the right side of (5.11), a formula that was proved,
whether or not there is a Siegel-Landau zero, by Friedlander [7]
using only sieve methods. Our proof of Proposition 5.2 can easily be
modified to work provided the quantity on the right side of (5.11)
is between $z/\phi(q)$ and $3z/\phi(q)$, for $z>q^B$ for some fixed
$B>0$; in fact our (non-elementary) proof gives such a result (as in
the discussion of the previous paragraph) but, better, Friedlander's
proof is easily modified to give this result, so that all of the
results in this paper can be proved avoiding use of zeros of
$L$-functions.
\endremark
\medskip

Lastly, we need a result which characterizes periodic functions that
are almost multiplicative.  We will show that such functions are
close to being characters. L. Babai, K. Friedl and A. Luk{\' a}cs
[1] have  explored such questions in greater generality recently,
but for the sake of completeness we provide a proof.

\proclaim{Proposition 5.3}  Let $0\le \epsilon < 1/2$.  Let $g:
({\Bbb Z}/q{\Bbb Z} )^* \to {\Bbb C}$ be a function with $g(1)= 1$,
and $|g(ab)-g(a)g(b)| \le \epsilon$ for all $a$ and $b$ coprime to
$q$.  Then there exists a character $\chi \pmod q$ such that
$|\chi(a) - g(a) | \le \epsilon/(1-2\epsilon)$, for all $(a,q)=1$.
\endproclaim
\demo{Proof} For any character $\chi \pmod q$, define
$$
{\hat g}(\chi)
= \sum_{ a\pmod q} g(a) \cbar(a).
$$
We claim that there exists a character $\chi \pmod q$ with $|{\hat g}(\chi)| \ge (1-2\epsilon)\phi(q)$.
Granting this claim, we see that for any $(a,q)=1$
$$
\align
|(g(a) -\chi(a)){\hat g}(\chi)| &=\Big| \sum_{b \pmod q } g(a)g(b)\cbar(b)-\chi(a)\sum_{ab \pmod q}
g(ab)\cbar(ab)\Big|  \\
&= \Big| \sum_{b\pmod q} (g(a)g(b)-g(ab)) \cbar(b) \Big| \le \epsilon \phi(q),\\
\endalign
$$
which proves the Proposition.

Note that
$$
\sum\Sb \chi \pmod q \endSb |{\hat g}(\chi)|^2 = \phi(q)
\sumstar_{a\pmod q} |g(a)|^2  \tag{5.12}
$$
(where, as usual, $\sumstar$ denotes a sum over integers coprime with $q$). Further,
$$
\align
\sum\Sb \chi\pmod q \endSb {\hat g}(\chi)^2 \overline{{\hat g}(\chi)}
&= \phi(q)
\sumstar\Sb a,b \endSb g(a) g(b) \overline{g(ab)} =
\phi(q) \sumstar_{a, b} |g(a)g(b)|^2 +E \\
\endalign
$$
where
$$
|E| \le \phi(q)\sumstar_{a,b} |g(ab)-g(a)g(b)||g(a)g(b)| \le
\epsilon \phi(q) \sumstar_{a,b} |g(a)g(b)| \le \epsilon \phi(q)^2
\sumstar_a |g(a)|^2,
$$
so that, by (5.12),
$$
\align \max_{\chi \pmod q} |{\hat g}(\chi)|  \cdot \phi(q) \sumstar_a
|g(a)|^2 &\geq \Big|\sum\Sb \chi\pmod q \endSb {\hat
g}(\chi)^2 \overline{{\hat g}(\chi)} \Big| \\
&\ge \phi(q) \Big(\sumstar_{a} |g(a)|^2\Big)^2 - \epsilon \phi(q)^2
\sumstar_a |g(a)|^2.
\endalign
$$
We deduce that there exists a $\chi \pmod q$ with $|g(\chi)| \ge
\sumstar_a |g(a)|^2 -\epsilon \phi(q)$. Now, using Cauchy's
inequality, we have
$$
\sumstar_{a} |g(a)|^2 \ge \Big| \sumstar_a g(a)g(a^{-1}) \Big| \ge
(1-\epsilon) \phi(q),
$$
since $|g(1) - g(a)g(a^{-1})| \le \epsilon$, so that we may conclude
the existence of a $\chi \pmod q$ with $|{\hat g}(\chi)|\ge
(1-2\epsilon) \phi(q)$, as claimed.
\enddemo

\demo{Proof of Theorem 1 for large moduli} Select $c \pmod q$ with
$(c,q)=1$ for which $|F(x;q,c) |$ is maximized. We may suppose that
$$
  |F(x;q,c) |\gg \frac{x}{q} \frac{1}{\sqrt{\log A}},
$$
else there is nothing to prove.  With this assumption, taking
$a=b=c$ in Proposition 5.2, and noting that
$|F(x;q,c^2) |\leq |F(x;q,c) |$, we deduce that
$$
|F(x;q,1)| \ge \frac 12  |F(x;q,c) |.\tag{5.13}
$$
Let us define $g$ by setting $F(x;q,a)=g(a) F(x;q,1)$. Then $g$ is
periodic with period $q$, $g(1)=1$, and by Proposition 5.2 we have
$$
|g(ab)-g(a)g(b)| = O\Big(\frac{1}{\sqrt{\log A}} \frac{x}{q
|F(x;q,1)|}\Big).
$$
By Proposition 5.3 it follows that there is a character $\chi \pmod
q$ such that for all $(a,q)=1$
$$
|g(a) - \chi(a)| = O\Big( \frac{1}{\sqrt{\log A}} \frac{x}{q|
F(x;q,1)|}\Big); \tag{5.14}
$$
that
is
$$
F(x;q,a)=\chi(a) F(x;q,1) +O\Big( \frac{x}{q} \ \frac{1}{\sqrt{\log
A}} \Big)
$$
which is (1.3e); and we also deduce that
$$
\align
\sum_{n\le x} f(n) \cbar(n) &= \sum_{a\pmod q} \cbar(a) F(x;q,a) =
F(x;q,1) \sum_{a\pmod q} \cbar(a) g(a) \\
&= \phi(q) F(x;q,1)+ O \Big( \frac{1}{\sqrt{\log A}} x\frac{\phi(q)}{q }\Big).\tag{5.15}\\
\endalign
$$

If $\chi\overline{\psi}$ is non-principal then
$$
{\Bbb D}_q(f,{ \chi}(n)n^{-it};x)^2 \ge \frac{1}{17} \log A,
$$
for all $|t|\le A$ by Lemma 3.4.   Therefore by Corollary  2.2 we
have
$$
\sum_{n\le x} f(n) \cbar(n) \ll  \frac{\phi(q)}q \ \frac{x}{A^{\frac
1{20}}} , \tag{5.16}
$$
which implies the result, using (5.13) and (5.15), when $r$ does not
divide $q$. If $r$  divides $q$, with $\chi\overline{\psi}$
non-principal, then, proceeding as above,
$$
\sum\Sb n\le x\\ (n,q)=1\endSb f(n) \overline{\psi(n)}= F(x;q,1)
\sum\Sb a\pmod q\\ (a,q)=1\endSb \overline{\psi(a)}g(a) =
O\Big(x\frac{\phi(q)}{q} \frac{1}{\sqrt{\log A}}\Big),
$$
by (5.14) and the orthogonality of the characters $\chi$ and $\psi$,
which gives the result.

If  $r|q$ and $\chi$ is induced by $\psi$ then from (5.15) we get
$$
\align F(x;q,a) &= g(a) F(x;q,1) = \chi(a) F(x;q,1) + O
\Big(\frac{x}{q} \frac{1}{\sqrt{\log A}}\Big)
\\
&= \frac{\psi(a)}{\phi(q)} \sum\Sb n\le x\\ (n,q)=1 \endSb f(n)
\overline{\psi(n)} + O\Big(\frac{x}{q} \frac{1}{\sqrt{\log A}}\Big),
\\
\endalign
$$
as desired.

\enddemo

 \enddocument